# Associated timelike helices in Minkowski 3-space

Mehmet Önder

*Delibekirli Village, Tepe Street, No: 63, 31440, Kırıkhan, Hatay, Turkey.*
*E-mail: mehmetonder197999@gmail.com*

**Abstract**
In this study, some new types of timelike general helices associated to a non-lightlike curve are introduced in Minkowski 3-space. These new helices are called associated timelike helices. Some special types of associated timelike helices are introduced and also by considering the conditions that reference curve is a non-lightlike helix or a spacelike slant helix, the position vectors of these new timelike helices are determinate.

**MSC:** 53A04, 53C40.

**Key words:** Non-lightlike helix; non-lightlike slant helix; position vector.

## 1. Introduction

Lorentzian geometry, closely connected with hyperbolic geometry, has an important place in geometry and physics. German mathematician Felix Klein is known for his works in non-Euclidean geometry comprising Lorentzian geometry. In 1873, the familiar scientist introduced Lorentzian geometry in his paper "Ueber die sogenannte Nicht-Euklidische Geometrie". Poincàre, known as mathematician and theoretical physicist, presented 3-dimensional Lorentzian geometry and described for the first time the Lorentzian transformations by the aid of their symmetrical form. In 1908, Lorentzian geometry was improved by Minkowski known for his works in relativity and Minkowski spaces[19].

From past to present, Lorentzian geometry has preserved its attraction and many scientists from different branches such as biology and mathematics have been dealing with the theory of curves related to Lorentzian geometry. From the viewpoint of differential geometry and the scientists, one of the most fascinating curves is helix and it has been studied a lot in different spaces by mathematicians and biologists interested in explanation of a DNA molecule [2,6,11,12,14,16,18,20,22,23].

A geometric non-lightlike curve is a general helix precisely when the ratio of its non-zero curvature to its non-zero torsion is a constant function. Much more important characterizations for a curve to be a general helix in the sense of Lorentzian have been given by Barros et al. and Ferrandez et al. in [9,14]. Furthermore, Ali has studied position vectors of general helices in Euclidean 3-space [4]. He has also given position vectors of spacelike general helices [2]. Later, Ali and Turgut have introduced position vectors of timelike general helices [6].

Moreover, Izumiya and Takeuchi have defined a new type of special curves called slant helix in Euclidean 3-space $E^3$ [15]. Ali has given position vectors of slant helices in $E^3$ [5]. Corresponding definitions and characterizations for slant helix in the sense of Lorentzian have been presented by Ali and Lopez [1]. Such a Lorentzian slant helix is identified by the fact that its principal normal vectors make a constant angle with a fixed line named its axis. A non-lightlike slant helix in $E_1^3$ has been characterized by the differential equation of its curvature $\kappa$ and its torsion $\tau$ given by

$$\frac{\kappa^2}{\left(\tau^2-\kappa^2\right)^{3/2}}\left(\frac{\tau}{\kappa}\right)' = \text{constant} \quad \text{or} \quad \frac{\kappa^2}{\left(\kappa^2\pm\tau^2\right)^{3/2}}\left(\frac{\tau}{\kappa}\right)' = \text{constant}$$

(see [1] for details).

In this work, our goal is to study a new type of non-lightlike helix associated to a non-lightlike curve. This new helix is named non-lightlike associated helix. We classify these curves according to the reference curve. If the main curve is a non-lightlike helix (or slant helix, respectively), then associated helix is called helix-connected (or slant helix-connected, respectively) non-lightlike associated helix. In the present paper especially, we study associated timelike helices. Additionally, we determine position vectors of associated timelike helices for some special cases.

## 2. Preliminaries

This section includes a brief summary of space curves and definitions of general helix, slant helix and Darboux helix in Minkowski 3-space $E_1^3$.

The real vector space $R^3$ endowed with standard flat Lorentzian metric defined by $\langle,\rangle = -dx_1^2 + dx_2^2 + dx_3^2$ is called three-dimensional Minkowski space and denoted by $E_1^3$, where $(x_1, x_2, x_3)$ is a rectangular coordinate system of $E_1^3$. The Lorentzian cross product (or vector product) of two vectors $a = (a_1, a_2, a_3)$, $b = (b_1, b_2, b_3)$ in $E_1^3$ is given by

$$a \times b = - \begin{vmatrix} -i & j & k \\ a_1 & a_2 & a_3 \\ b_1 & b_2 & b_3 \end{vmatrix}.$$

A vector $a \in E_1^3$ is called a spacelike vector if $\langle a,a \rangle > 0$ or $a = 0$; called a timelike vector if $\langle a,a \rangle < 0$ and called a null(lightlike) vector if $\langle a,a \rangle = 0, a \neq 0$. The norm of a vector $a \in E_1^3$ is given by $\|a\| = \sqrt{|\langle a,a \rangle|}$. Similarly, a curve $\alpha(t): I \subset \mathbb{R} \to E_1^3$ is called a spacelike, timelike or null(lightlike) curve if all of its velocity vectors $\alpha'(t)$ are spacelike, timelike or null(lightlike), respectively[17]. The curves in $E_1^3$ have relationships with relativity in physics. Timelike, spacelike and null curves correspond to the path of a particular moving at less than the speed of light, faster than the speed of light and at a speed of light, respectively [13].

Let $\{T, N, B\}$ denotes the Frenet frame of a non-lightlike curve $\alpha$. Then, the curve $\alpha$ is a timelike or a spacelike curve. The spacelike curves with non-lightlike Frenet vectors have two types according to the Lorentzian casual characters of Frenet vectors: A spacelike curve is called of type 1(respectively, type 2) if its principal normal vector $N$ (respectively, binormal vector $B$) is timelike and other Frenet vectors are spacelike [10]. For the derivatives of Frenet frame of a unit speed non-lightlike curve $\alpha(s)$, the Frenet-Serret formulas are given separately in [21] according to Lorentzian casual characters of space curves. We can give these formulas with a single equality for each case as follows:

$$\begin{bmatrix} T' \\ N' \\ B' \end{bmatrix} = \begin{bmatrix} 0 & \kappa & 0 \\ \varepsilon_B \kappa & 0 & \tau \\ 0 & \varepsilon_T \tau & 0 \end{bmatrix} \begin{bmatrix} T \\ N \\ B \end{bmatrix}, \qquad (1)$$

where $\varepsilon_Y = \langle Y, Y \rangle$, $\varepsilon_Y^2 = 1$, $Y \in \{T, N, B\}$, $\kappa(s)$ is the curvature(or first curvature) and $\tau(s)$ is the torsion (or second curvature) of $\alpha$ at $s$ defined by $\tau(s) = \varepsilon_N \varepsilon_B \langle B', N \rangle$. The Frenet vectors satisfy the relations

$$B = \varepsilon_T \varepsilon_N T \times N, \quad N = \varepsilon_B \varepsilon_T B \times T, \quad T = \varepsilon_N \varepsilon_B N \times B, \quad \varepsilon_B = -\varepsilon_T \varepsilon_N.$$

The vector $W$ defined by $W(s) = -\varepsilon_B \tau(s) T(s) - \varepsilon_N \kappa(s) B(s)$ is called Darboux vector of non-lightlike curve $\alpha$. This vector allows us to write Frenet formulas (1) in a different form such as $T' = W \times T$, $N' = W \times N$, $B' = W \times B$.

Unlike Euclidean space, in Minkowski 3-space, four different angle between two non-lightlike vectors have been given in [19] as follows:

**Definition 2.1.** Let $x$ and $y$ be spacelike vectors in $E_1^3$ that span a spacelike vector subspace. Then, we have $|\langle x, y \rangle| \leq \|x\| \|y\|$ and hence, there is a unique real number $\theta \geq 0$ such that $|\langle \vec{x}, \vec{y} \rangle| = \|\vec{x}\| \|\vec{y}\| \cos\theta$. This number is called the Lorentzian spacelike angle between the vectors $x$ and $y$.

**Definition 2.2.** Let $x$ and $y$ be spacelike vectors in $E_1^3$ that span a timelike vector subspace. Then, we have $|\langle x, y \rangle| > \|x\| \|y\|$ and hence, there is a unique real number $\theta \geq 0$ such that $|\langle x, y \rangle| = \|x\| \|y\| \cosh\theta$. This number is called the Lorentzian timelike angle between the vectors $x$ and $y$.

**Definition 2.3.** Let $x$ be a spacelike vector and $y$ be a positive timelike vector in $E_1^3$. Then, there is a unique non-negative real number $\theta$ such that $|\langle x, y \rangle| = \|x\| \|y\| \sinh\theta$. This number is called the Lorentzian timelike angle between the vectors $x$ and $y$.

**Definition 2.4.** Let $x$ and $y$ be positive(negative) timelike vectors in $E_1^3$. Then, there is a unique non-negative real number $\theta$ such that $\langle x, y \rangle = \|x\| \|y\| \cosh\theta$. This number is called the Lorentzian timelike angle between the vectors $x$ and $y$.

Let $\psi = \psi(s)$ be a naturel parametrization of a non-lightlike unit speed regular curve in $E_1^3$. Then, $\psi$ is called an $i$-slant helix if the unit vector
$$\psi_{i+1} = \frac{\psi_i'(s)}{\|\psi_i'(s)\|},$$
makes a constant angle with a constant vector, where $\psi_0 = \psi(s)$ [8]. Moreover, in [24], the authors defined Darboux helix in Euclidean 3-space. We give the corresponding definition for non-lightlike curves as follows: A non-lightlike curve $\alpha$ is called a non-lightlike Darboux helix if the Darboux vector $W$ of curve always makes a constant angle with a fixed non-lightlike direction in Minkowski 3-space.

**3. Associated timelike helices in Minkowski 3-space**

This section includes the definitions of associated timelike helices and classifications of those curves for some special cases. Also, the position vectors of associated timelike helices are determinate.

Let $\alpha : I \to E_1^3$ be a unit speed non-lightlike curve with arc-length parameter $s$ and let $\{T, N, B\}$ and $\kappa(s)$, $\tau(s)$ denote Frenet frame and non-zero curvature functions of $\alpha$, respectively. Consider non-lightlike curve $\beta : J \to E_1^3$ given by the parametrization
$$\beta(s) = \alpha(s) + a_1(s) T(s) + a_2(s) N(s) + a_3(s) B(s), \tag{2}$$

where $a_i = a_i(s)$; $(1 \leq i \leq 3)$ are differentiable functions of arc-length parameter $s$. Then, the curve $\beta$ is called a non-lightlike associated curve of $\alpha$. This definition is a generalization for non-lightlike curve pairs and considering with the special properties, parametric forms of well-known curve pairs can be given by (2). For instance, if we choose $a_1 = a_3 = 0$, $a_2 = \lambda$ is a non-zero constant, and assume that the unit principal normals of $\alpha$ and $\beta$ are linearly dependent, then we obtain the parametric form of a Bertrand mate of non-lightlike curve $\alpha$. Similarly, with proper chosen of functions $a_i = a_i(s)$; $(1 \leq i \leq 3)$, parametric forms of other famous curve pairs such us Mannheim-partner curves and involute-evolute curves can be given by (2).

If $\beta$ in (2) is a general helix, then it is called a non-lightlike associated helix connected to $\alpha$. The function $d(s)$ defined by $d(s) = \|\beta(s) - \alpha(s)\|$ is called distance function between the corresponding points of non-lightlike associated curves $\alpha$ and $\beta$.

If we differentiate (2) with respect to $s$ and consider (1), we have
$$\beta'(s) = F_1(s)T(s) + F_2(s)N(s) + F_3(s)B(s), \tag{3}$$
where $F_i = F_i(s)$; $(1 \leq i \leq 3)$ are differentiable functions of $s$ defined by
$$F_1 = 1 + a_1' + \varepsilon_B a_2 \kappa, \quad F_2 = a_2' + a_1 \kappa + \varepsilon_T a_3 \tau, \quad F_3 = a_3' + a_2 \tau. \tag{4}$$
Equations (3) and (4) allow us to introduce and investigate non-lightlike special associated helices in Minkowski 3-space. Privately, in this study we study associated timelike helices.

**3.1. Helix-Connected Associated(HCA) Timelike Helices**

Let $\alpha$ be a non-lightlike general helix and $\beta$ is a non-lightlike associated curve of $\alpha$. Then, $\beta$ is called helix-connected non-lightlike associated curve or HC-NL-associated curve. Since we assume $\alpha$ is a non-lightlike general helix, we can write $\tau(s) = m\kappa(s)$, where $m$ is a non-zero constant. Let now consider special case that $\beta$ is a timelike curve such that tangent vector $\beta'$ of $\beta$ is linearly dependent with unit tangent vector $T$ of $\alpha$. Then, we have that $\alpha$ is also a timelike curve. For this case, we have $F_1 \neq 0$, $F_2 = F_3 = 0$ and from (4), the following system is obtained
$$1 + a_1' + a_2 \kappa \neq 0, \quad a_2' + a_1 \kappa - a_3 \tau = 0, \quad a_3' + a_2 \tau = 0. \tag{5}$$
Now, the following special cases can be given:

**Case (1):** $a_1 = 0$. In this case, system (5) becomes,
$$1 + a_2 \kappa \neq 0, \quad a_2' - a_3 \tau = 0, \quad a_3' + a_2 \tau = 0, \tag{6}$$
and solving these equations gives
$$a_1 = 0, \quad a_2 = \sin\left(\int \tau(s)ds\right), \quad a_3 = \cos\left(\int \tau(s)ds\right). \tag{7}$$
By writing (7) in (2), the position vector of associated timelike curve $\beta$ is obtained as
$$\beta(s) = \alpha(s) + \sin\left(\int \tau(s)ds\right)N(s) + \cos\left(\int \tau(s)ds\right)B(s), \tag{8}$$
which gives the followings:

**Theorem 3.1.** *The associated curve $\beta$ given in (8) is a timelike helix if and only if $\alpha$ is a timelike helix.*

**Remark 3.2.** *The associated curve (8) can be named by: Helix-connected associated timelike helix of type 1 or HCA-timelike helix of type 1.*

**Case (2):** $a_2 = 0$. In this case, from (5) we have $a_1 = c\dfrac{\tau}{\kappa}$ and $a_3 = c$, where $c$ is constant. If we take $c = 0$, it follows $a_i = 0$; $(1 \leq i \leq 3)$ which gives $\beta(s) = \alpha(s)$. So, we take $c$ is non-zero and from (2), the position vector of associated curve $\beta$ is obtained as

$$\beta(s) = \alpha(s) + c\left(\dfrac{\tau(s)}{\kappa(s)}T(s) + B(s)\right), \tag{9}$$

and the followings are obtained:

**Theorem 3.3.** *The associated curve $\beta$ given in (9) is a timelike helix if and only if $\alpha$ is a timelike helix.*

**Remark 3.4.** *The associated curve (9) can be named by: Helix-connected associated timelike helix of type 2 or HCA-timelike helix of type 2.*

**Case (3):** $a_3 = 0$. In this case, from system (5) we have $a_1 = a_2 = 0$ and from (2), it follows $\beta = \alpha$.

On the other hand, we know that if $\alpha$ is a spacelike general helix, its binormal vector $B(s)$ makes a constant angle with a constant vector. For this situation, assume that $\alpha$ is a spacelike curve of type 2. Then, unit binormal vector $B(s)$ is timelike and we can consider the situation that the tangent of associated curve $\beta$ is linearly dependent with timelike binormal vector $B(s)$ of spacelike curve $\alpha$. For this case, from (4), it follows

$$a_3' + a_2\tau \neq 0,\ 1 + a_1' - a_2\kappa = 0,\ a_2' + a_1\kappa + a_3\tau = 0, \tag{10}$$

and following special cases can be written:

**Case (4):** $a_1 = 0$. Then, system (10) reduces to

$$a_3' + a_2\tau \neq 0,\ 1 - a_2\kappa = 0,\ a_2' + a_3\tau = 0,$$

and the solution is

$$a_1 = 0,\ a_2(s) = \dfrac{1}{\kappa(s)},\ a_3(s) = -\dfrac{1}{\tau(s)}\left(\dfrac{1}{\kappa(s)}\right)'.$$

Then, from (2), the parametric form of associated curve $\beta$ is

$$\beta(s) = \alpha(s) + \dfrac{1}{\kappa(s)}N(s) - \dfrac{1}{\tau(s)}\left(\dfrac{1}{\kappa(s)}\right)' B(s), \tag{11}$$

and followings can be given:

**Theorem 3.5.** *The associated curve $\beta$ given in (11) is a timelike general helix if and only if $\alpha$ is a spacelike general helix with unit timelike binormal vector $B(s)$.*

**Remark 3.6.** *The associated curve (11) can be named by: Helix-connected associated timelike helix of type 3 or HCA-timelike helix of type 3.*

**Case (5):** $a_2 = 0$. In this case, from (10) we have,

$$a_3' \neq 0,\ 1 + a_1' = 0,\ a_1\kappa + a_3\tau = 0$$

and the solution is

$$a_1(s) = \upsilon - s,\ a_2 = 0,\ a_3(s) = -(\upsilon - s)\dfrac{\kappa(s)}{\tau(s)}.$$

where $\upsilon$ is integration constant. Then, from (2), the parametric form of associated curve $\beta$ is

$$\beta(s) = \alpha(s) + (\upsilon - s)\left(T(s) - \frac{\kappa(s)}{\tau(s)} B(s)\right). \tag{12}$$

Now, we can write followings:

**Theorem 3.7.** *The associated curve $\beta$ given in (12) is a timelike general helix if and only if $\alpha$ is a spacelike general helix with unit timelike binormal vector $B(s)$.*

**Remark 3.8.** *The associated curve (12) can be named by: Helix-connected associated helix of type 4 or HCA-timelike helix of type 4.*

**Case (6) :** $a_3 = 0$. For this case, system (10) gives us

$$a_2 \tau \neq 0, \ 1 + a_1' - a_2 \kappa = 0, \ a_2' + a_1 \kappa = 0. \tag{13}$$

From the last equation, we obtain

$$a_1' = -\frac{1}{\kappa} a_2'' - \left(\frac{1}{\kappa}\right)' a_2', \tag{14}$$

and by writing this result in the second equation in (13), it follows

$$a_2'' + \kappa \left(\frac{1}{\kappa}\right)' a_2' + \kappa^2 a_2 - \kappa = 0. \tag{15}$$

By changing the variable as $\mu = \int \kappa(s) ds$, the homogeny part of (15) becomes

$$\frac{d^2 a_2}{d\mu^2} + a_2 = 0,$$

and it has the solution

$$a_{2h}(s) = c_1 \cos\left(\int \kappa(s) ds\right) + c_2 \sin\left(\int \kappa(s) ds\right), \tag{16}$$

where $c_1, c_2$ are real constants. Now, considering (16) and by using the variation of parameters method, the solution of (15) is

$$\begin{aligned} a_2(s) = &c_1 \cos\left(\int \kappa(s) ds\right) + c_2 \sin\left(\int \kappa(s) ds\right) - \cos\left(\int \kappa(s) ds\right) \int \left[\sin\left(\int \kappa(s) ds\right)\right] ds \\ &+ \sin\left(\int \kappa(s) ds\right) \int \left[\cos\left(\int \kappa(s) ds\right)\right] ds \end{aligned} \tag{17}$$

Differentiating (17) and writing the result in the second equation in (13), we have

$$\begin{aligned} a_1(s) = &-\sin\left(\int \kappa(s) ds\right) \left[\int \left[\sin\left(\int \kappa(s) ds\right)\right] ds - c_1\right] \\ &- \cos\left(\int \kappa(s) ds\right) \left[\int \left[\cos\left(\int \kappa(s) ds\right)\right] ds + c_2\right] \end{aligned} \tag{18}$$

Then, the parametric form of associated curve $\beta$ is determinate by

$$\beta(s) = \alpha(s) + a_1(s) T(s) + a_2(s) N(s), \tag{19}$$

where $a_1(s)$ and $a_2(s)$ are given in (18) and (17), respectively. Now, we have the followings:

**Theorem 3.9.** *The associated curve $\beta$ given in (19) is a timelike general helix if $\alpha$ is a spacelike general helix with unit timelike binormal vector $B(s)$.*

**Remark 3.10.** *The associated curve (19) can be named by: Helix-connected associated timelike helix of type 5 or HCA-timelike helix of type 5.*

Furthermore, let $d_i(s)$ denotes the distance function for HCA-timelike helix of type $i$. Then, for $d_5(s)$, we have $d_5(s) = \sqrt{a_1^2(s) + a_2^2(s)}$. It is clear that $d_5(s)$ is constant if and only if $a_1 a_1' + a_2 a_2' = 0$. On the other hand, multiplying the second and third equalities in (13) with $a_1$ and $a_2$, respectively, and adding the results gives $a_1 a_1' + a_2 a_2' = -a_1$. Then, it follows $d_5(s)$ is constant if and only if $a_1 = 0$, or equivalently $a_2$ is constant. Similarly, by considering the differentiation of other distance functions, we can give the followings:

**Corollary 3.11.** *(i) The distance functions $d_1$ and $d_2$ are constants at all points.*
*(ii) The distance function $d_3$ is constant if and only if $\kappa(s)$ is constant or*
$$\frac{1}{\kappa(s)}\left(\frac{1}{\kappa(s)}\right)' = m^2 s + c,$$
*where $c$ is integration constant.*
*(iii) The distance function $d_4$ is constant and equal to zero at the intersection points $\beta(\upsilon) = \alpha(\upsilon)$ and non-constant all other points.*
*(iv) The distance function $d_5$ is constant if and only if $a_1 = 0$ or $a_2$ is constant.*

Let now determinate the position vector of a HCA-timelike helix. For this, first assume that main curve $\alpha$ is a timelike helix with $|\tau(s)/\kappa(s)| < 1$. Without loss of generality, we can choose the spacelike base vector $e_3 = (0,0,1)$ as the axis of timelike helix $\alpha$. Then, we have the following theorem:

**Theorem 3.12.**([6]). *The position vector $\psi$ of a timelike general helix whose tangent vector makes a constant Lorentzian timelike angle with a fixed spacelike straight line in the space, is computed in the natural representation form:*
$$\psi(s) = \sqrt{1+n^2} \int \left( \cosh\left[\sqrt{1-m^2}\int \kappa(s)ds\right], \sinh\left[\sqrt{1-m^2}\int \kappa(s)ds\right], m \right) ds \qquad (20)$$
*or in the parametric form as follows:*
$$\psi(\theta) = \int \frac{\sqrt{1+n^2}}{\kappa(\theta)} \left( \cosh\left[\sqrt{1-m^2}\,\theta\right], \sinh\left[\sqrt{1-m^2}\,\theta\right], m \right) d\theta \qquad (21)$$
*where $\theta = \int \kappa(s)ds$, $m = n/\sqrt{1+n^2}$, $n = \sinh(\varphi)$ and $\varphi$ is the timelike angle between the fixed spacelike straight line $e_3$ (axis of timelike general helix) and the tangent vector of the curve $\psi$.*

From (20), Frenet vectors of timelike helix $\alpha$ are calculated as follows:
$$T(s) = \sqrt{1+n^2}\left( \cosh\left[\sqrt{1-m^2}\int \kappa(s)ds\right], \sinh\left[\sqrt{1-m^2}\int \kappa(s)ds\right], m \right),$$
$$N(s) = \left( \sinh\left[\sqrt{1-m^2}\int \kappa(s)ds\right], \cosh\left[\sqrt{1-m^2}\int \kappa(s)ds\right], 0 \right),$$
$$B(s) = \sqrt{1+n^2}\left( m\cosh\left[\sqrt{1-m^2}\int \kappa(s)ds\right], m\sinh\left[\sqrt{1-m^2}\int \kappa(s)ds\right], 1 \right).$$

Now, we can give the followings for the special cases 1 and 2:

**Theorem 3.13.** *The position vector $\beta(s) = (\beta_1, \beta_2, \beta_3)$ of the HCA-timelike helix of type 1 with spacelike axis $e_3$, is obtained as*

$$\beta_1(s) = \sqrt{1+n^2} \int \cosh\left[\sqrt{1-m^2}\int \kappa(s)ds\right]ds + \sin\left[m\int \kappa(s)ds\right]\sinh\left[\sqrt{1-m^2}\int \kappa(s)ds\right]$$
$$+ n\cos\left[m\int \kappa(s)ds\right]\cosh\left[\sqrt{1-m^2}\int \kappa(s)ds\right],$$

$$\beta_2(s) = \sqrt{1+n^2} \int \sinh\left[\sqrt{1-m^2}\int \kappa(s)ds\right]ds + \sin\left[m\int \kappa(s)ds\right]\cosh\left[\sqrt{1-m^2}\int \kappa(s)ds\right]$$
$$+ n\cos\left[m\int \kappa(s)ds\right]\sinh\left[\sqrt{1-m^2}\int \kappa(s)ds\right],$$

$$\beta_3(s) = (ns+l) + \sqrt{1+n^2}\cos\left[m\int \kappa(s)ds\right],$$

*where $l$ is integration constant.*

**Theorem 3.14.** *The position vector $\beta(s) = (\beta_1, \beta_2, \beta_3)$ of the HCA-timelike helix of type 2 with spacelike axis $e_3$, is obtained as*

$$\beta_1(s) = \sqrt{1+n^2}\int \cosh\left[\sqrt{1-m^2}\int \kappa(s)ds\right]ds,$$
$$\beta_2(s) = \sqrt{1+n^2}\int \sinh\left[\sqrt{1-m^2}\int \kappa(s)ds\right]ds,$$
$$\beta_3(s) = (ns+l) + c\left(nm - \sqrt{1+n^2}\right),$$

*where $l$ is integration constant.*

For $n = \sqrt{3}/3$, $m = 1/2$, the graphs of timelike general helix $\alpha$ with curvatures $\kappa = 6$, $\tau = 3$ and HCA-timelike helices of type 1 and of type 2 with $c = 1$ are given in left side, middle and right sight of Fig. 1, respectively.

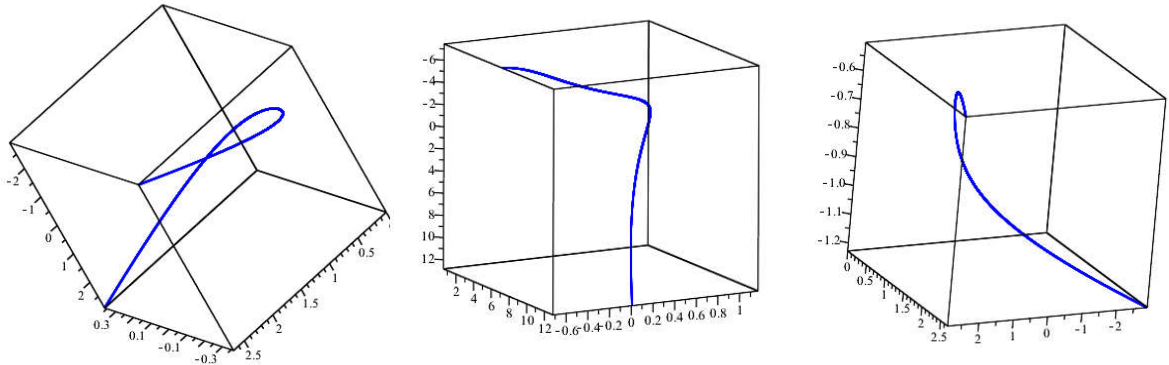

**Fig. 1.** From left to right, helix $\alpha$ and HCA-timelike helices of type 1 and 2, respectively.

In Theorem 3.12, if the curve $\psi$ is a timelike helix with a timelike axis $e_1$, then the natural representation form of $\psi$ is

$$\psi(s) = \sqrt{n^2-1}\int \left(m, \cos\left[\sqrt{m^2-1}\int \kappa(s)ds\right], \sin\left[\sqrt{m^2-1}\int \kappa(s)ds\right]\right)ds,$$

where $\theta = \int \kappa(s)ds$, $m = n/\sqrt{n^2-1}$, $n = \cosh(\varphi)$ and $\varphi$ is the constant Lorentzian timelike angle between the fixed timelike straight line $e_1$ (axis of timelike general helix) and the tangent vector of the curve $\psi$ [6]. Then, corresponding theorems for HCA-timelike helices of type 1 and 2 can be given similarly.

If the curve $\psi$ is a spacelike helix with a timelike binormal vector, then there exists two situations [2]:

ii) If $\psi$ is a spacelike curve of type 2 with $|\tau(s)/\kappa(s)|>1$, the axis is spacelike and constant angle between the unit tangent of curve and axis is a timelike angle, then the natural representation form of $\psi$ is

$$\psi(s) = \sqrt{n^2-1} \int \left( \cosh\left[\sqrt{m^2-1}\int \kappa(s)ds\right], \sinh\left[\sqrt{m^2-1}\int \kappa(s)ds\right], m \right) ds,$$

where $\theta = \int \kappa(s)ds$, $m = n/\sqrt{n^2-1}$, $n = \cosh(\varphi)$ and $\varphi$ is the timelike angle between the fixed spacelike straight line $e_3$ (axis of spacelike general helix of type 1) and the tangent vector of the curve $\psi$.

iii) If $\psi$ is a spacelike curve of type 2 with $|\tau(s)/\kappa(s)|<1$, the axis is timelike and constant angle between the unit tangent of curve and axis is a timelike angle, then the natural representation form of $\psi$ is

$$\psi(s) = \sqrt{n^2+1} \int \left( m, \cos\left[\sqrt{1-m^2}\int \kappa(s)ds\right], \sin\left[\sqrt{1-m^2}\int \kappa(s)ds\right] \right) ds,$$

where $\theta = \int \kappa(s)ds$, $m = n/\sqrt{1+n^2}$, $n = \sinh(\varphi)$ and $\varphi$ is the timelike angle between the fixed timelike straight line $e_1$ (axis of spacelike general helix of type 2) and the tangent vector of the curve $\psi$. Then, from (11), (12) and (19), the position vectors of HCA-timelike helices of types 3, 4 and 5 can be determinate similarly.

**3.2. Slant helix-connected associated (SHCA) timelike helices**

Assume that the tangent vector $\beta'$ of associated curve $\beta$ is linearly dependent with unit principal normal vector $N$ of main curve $\alpha$. Then for Lorentzian causal character of $\beta$ we have followings:

i) $\beta$ is a timelike curve if and only if $\alpha$ is a spacelike curve of type 1.

ii) $\beta$ is a spacelike curve if and only if $\alpha$ is a timelike curve or a spacelike curve of type 2.

Now, from (3), we can write $\beta'(s) = F_2(s)N(s)$. Then, for the arc-length parameter $s_\beta$ of curve $\beta$, we have $ds_\beta = \pm F_2(s)ds$. After some computations, Frenet vectors of associated curve $\beta$ are calculated as follows:

$$\begin{cases} T_\beta(s) = \pm N(s), \ N_\beta(s) = \dfrac{1}{\sqrt{|\varepsilon_T \kappa^2(s)+\varepsilon_B \tau^2(s)|}}\left(\varepsilon_B \kappa(s)T(s)+\tau(s)B(s)\right), \\ B_\beta(s) = \dfrac{\mp 1}{\sqrt{|\varepsilon_T \kappa^2(s)+\varepsilon_B \tau^2(s)|}}\left(-\varepsilon_B \tau(s)T(s)-\varepsilon_N \kappa(s)B(s)\right) = \dfrac{\mp W(s)}{\sqrt{|\varepsilon_T \kappa^2(s)+\varepsilon_B \tau^2(s)|}}, \end{cases} \quad (22)$$

where $W(s)$ is Darboux vector of main curve $\alpha$. Now, the following theorem can be given:

**Theorem 3.15.** *Let $\beta$ be a non-lightlike associated curve such that tangent vector $\beta'$ of $\beta$ is linearly dependent with unit principal normal vector $N$ of main curve $\alpha$. The followings are equivalent,*

*(i)* $\beta$ *is a non-lightlike helix.*
*(ii)* $\alpha$ *is a non-lightlike slant helix.*
*(iii)* $\alpha$ *is a non-lightlike Darboux helix.*

**Remark 3.16.** *The helix $\beta$ associated to slant helix $\alpha$ can be named by: Slant helix-connected associated non-lightlike helix or SHCA-non-lightlike helix.*

Now, we will consider the case that $\beta$ is a timelike curve and introduce special types of SHCA-timelike helices. For this purpose, let assume that $\alpha$ is a spacelike curve of type 1, i.e., principal normal $N(s)$ is timelike. Then, from (3) we have $F_2 \neq 0$, $F_1 = F_3 = 0$ and from (4), the following system is obtained

$$a_2' + a_1\kappa + a_3\tau \neq 0, \ 1 + a_1' + a_2\kappa = 0, \ a_3' + a_2\tau = 0. \tag{23}$$

Now, we can consider the following special cases:

**Case (A):** $a_1 = 0$. In this case, system (23) becomes

$$a_2' + a_3\tau \neq 0, \ 1 + a_2\kappa = 0, \ a_3' + a_2\tau = 0. \tag{24}$$

Solving these equations, it follows

$$a_1 = 0, \ a_2 = -\frac{1}{\kappa(s)}, \ a_3 = \int \frac{\tau(s)}{\kappa(s)} ds. \tag{25}$$

By writing (25) in (2), the position vector of associated timelike curve $\beta$ is obtained as

$$\beta(s) = \alpha(s) - \frac{1}{\kappa(s)} N(s) + \left[ \int \frac{\tau(s)}{\kappa(s)} ds \right] B(s), \tag{26}$$

and following results are obtained:

**Theorem 3.17.** *The associated curve $\beta$ given in (26) is a timelike general helix if and only if $\alpha$ is a spacelike slant helix with timelike principal normal $N(s)$.*

**Remark 3.18.** *The associated curve (26) can be named by: Slant helix-connected associated timelike helix of type 1 or SHCA-timelike helix of type 1.*

**Case (B):** $a_2 = 0$. Then, system (23) becomes

$$a_1\kappa + a_3\tau \neq 0, \ 1 + a_1' = 0, \ a_3' = 0, \tag{27}$$

and it has the solution

$$a_1 = \xi - s, \ a_2 = 0, \ a_3 = \zeta. \tag{28}$$

where $\xi, \zeta$ are real constants. By writing (28) in (2), the position vector of timelike associated curve $\beta$ is obtained as

$$\beta(s) = \alpha(s) + (\xi - s)T(s) + \zeta B(s), \tag{29}$$

and following theorem is obtained:

**Theorem 3.19.** *The associated curve $\beta$ given in (29) is a timelike helix if and only if $\alpha$ is a spacelike slant helix with timelike principal normal $N(s)$.*

**Remark 3.20.** *The associated curve (29) can be named by: Slant helix-connected associated timelike helix of type 2 or SHCA-timelike helix of type 2.*

**Case (C):** $a_3 = 0$. From (23), it follows

$$a_2' + a_1\kappa \neq 0, \ 1 + a_1' + a_2\kappa = 0, \ a_2\tau = 0, \tag{30}$$

which has the solution

$$a_1 = \omega - s, \quad a_2 = 0, \quad a_3 = 0. \qquad (31)$$

where $\omega$ is integration constant. Writing (31) in (2) gives the position vector of associated curve $\beta$ as follows,

$$\beta(s) = \alpha(s) + (\omega - s)T(s). \qquad (32)$$

Then, the following theorem is given:

**Theorem 3.21.** *The associated curve $\beta$ given in (32) is a timelike helix if and only if $\alpha$ is a spacelike slant helix with timelike principal normal $N(s)$.*

**Remark 3.22.** *The associated curve (32) can be named by: Slant helix-connected associated timelike helix of type 3 or SHCA-timelike helix of type 3.*

Furthermore, by considering (26), (29) and (32), the followings can be given:

**Corollary 3.23.** *SHCA-timelike helix of type 2 is also a SHCA-timelike helix of type 3 if and only if $\zeta = 0$.*

Let $d_1$, $d_2$ and $d_3$ denote the distance functions for timelike associated curves (26), (29) and (32), respectively. Then, we have

$$d_1(s) = \sqrt{\left| -\left(\frac{1}{\kappa(s)}\right)^2 + \left[\int \frac{\tau(s)}{\kappa(s)} ds\right]^2 \right|}, \quad d_2(s) = \sqrt{(\xi - s)^2 + \zeta^2},$$

$$d_3(s) = \sqrt{|(\omega - s)^2|}.$$

By differentiating these functions, the following corollary can be given:

**Corollary 3.24.** *(i) For the distance function $d_1(s)$, the followings are equivalent:*

*(a) $d_1 = c$ is constant. (b) $\tau \left[\int \frac{\tau(s)}{\kappa(s)} ds\right] - \left(\frac{1}{\kappa}\right)' = 0$. (c) $\tau = \frac{g'}{\sqrt{c^2 + g^2}}$, where $g(s) = \frac{1}{\kappa(s)}$.*

*(ii) $d_2$ is constant and equal to $\zeta$ at points $\alpha(\xi)$ and non-constant at all other points.*

*(iii) $d_3$ is constant and equal to zero at the intersection points $\beta(\omega) = \alpha(\omega)$ and non-constant at all other points.*

Let now determine the position vector of a SHCA-timelike helix. For this, first assume that main curve $\alpha$ is a spacelike slant helix with timelike principal normal $N(s)$. Without loss of generality, we can choose the spacelike base vector $e_3 = (0, 0, 1)$ as the axis of spacelike slant helix $\alpha$. Then, we have the following theorem:

**Theorem 3.25.**([7]). *The position vector $\psi(s)$ of a spacelike slant helix, whose timelike principal normal vector makes a constant Lorentzian timelike angle with a fixed spacelike straight line $e_3$ is computed in the natural representation form:*

$$\psi(s) = \left(\frac{n}{m}\int\left[\int \kappa(s)\cosh[A]ds\right]ds, \frac{n}{m}\int\left[\int \kappa(s)\sinh[A]ds\right]ds, n\int\left[\int \kappa(s)ds\right]ds\right) \qquad (33)$$

where $A = \frac{1}{n}\arcsin\left(m\int \kappa(s)ds\right)$, $\tau(s) = \pm \frac{m\kappa(s)\int \kappa(s)ds}{\sqrt{1-m^2\left(\int \kappa(s)ds\right)^2}}$, $\theta = \int \kappa(s)ds$, $m = \frac{n}{\sqrt{1+n^2}}$,

$n = \sinh(\varphi)$ and $\varphi$ is the constant timelike angle between the axis of spacelike slant helix and timelike principal normal vector of the curve $\psi$.

From (33) the Frenet vectors of slant helix $\alpha$ are computed as follows:

$$T(s) = \left(\frac{n}{m}\int \kappa(s)\cosh[A]ds, \frac{n}{m}\int \kappa(s)\sinh[A]ds, n\int \kappa(s)ds\right) \quad (34)$$

$$N(s) = \left(\frac{n}{m}\cosh[A], \frac{n}{m}\sinh[A], n\right), \quad (35)$$

$$B(s) = \left(\frac{n^2}{m}\sinh[A]\int \kappa(s)ds - \frac{n^2}{m}\int \kappa(s)\sinh[A]ds,\right.$$

$$+\frac{n^2}{m}\cosh[A]\int \kappa(s)ds - \frac{n^2}{m}\int \kappa(s)\cosh[A]ds, \quad (36)$$

$$\left.+\frac{n^2}{m^2}\sinh[A]\int \kappa(s)\cosh[A]ds - \frac{n^2}{m^2}\cosh[A]\int \kappa(s)\sinh[A]ds\right)$$

Now, by writing (34)-(36) in (26), (29) and (32), respectively, we can give the following theorems:

**Theorem 3.26.** *The position vector $\beta(s) = (\beta_1, \beta_2, \beta_3)$ of the SHCA-timelike helix of type 1 with spacelike axis $e_3$ and constant spacelike angle $\varphi$, between the unit tangent of SHCA-timelike helix and $e_3$, is obtained as*

$$\beta_1(s) = \frac{n}{m}\int\left[\int \kappa(s)\cosh[A]ds\right]ds - \frac{n}{m\kappa(s)}\cosh[A]$$

$$-\frac{n^2}{m}(ms+c_1)\left(\int \kappa(s)\sinh[A]ds - \sinh[A]\int \kappa(s)ds\right),$$

$$\beta_2(s) = \frac{n}{m}\int\left[\int \kappa(s)\sinh[A]ds\right]ds - \frac{n}{m\kappa(s)}\sinh[A]$$

$$-\frac{n^2}{m}(ms+c_2)\left(\int \kappa(s)\cosh[A]ds - \cosh[A]\int \kappa(s)ds\right),$$

$$\beta_3(s) = n\int\left[\int \kappa(s)ds\right]ds - \frac{n}{\kappa(s)}$$

$$+\frac{n^2}{m^2}(ms+c_3)\left(\sinh[A]\int \kappa(s)\cosh[A]ds - \cosh[A]\int \kappa(s)\sinh[A]ds\right),$$

*where $c_i$, $(1 \leq i \leq 3)$ are integration constants.*

**Theorem 3.27.** *The position vector $\beta(s) = (\beta_1, \beta_2, \beta_3)$ of the SHCA-timelike helix of type 2 with spacelike axis $e_3$ and constant spacelike angle $\varphi$, between the unit tangent of SHCA-timelike helix and $e_3$, is obtained as*

$$\beta_1(s) = \frac{n}{m}\int\left[\int \kappa(s)\cosh[A]ds\right]ds + \frac{n}{m}(\xi - s)\int \kappa(s)\cosh[A]ds$$
$$+ \zeta\frac{n^2}{m}\left(\sinh[A]\int \kappa(s)ds - \int \kappa(s)\sinh[A]ds\right),$$
$$\beta_2(s) = \frac{n}{m}\int\left[\int \kappa(s)\sinh[A]ds\right]ds + \frac{n}{m}(\xi - s)\int \kappa(s)\sinh[A]ds$$
$$+ \zeta\frac{n^2}{m}\left(\cosh[A]\int \kappa(s)ds - \int \kappa(s)\cosh[A]ds\right),$$
$$\beta_3(s) = n\int\left[\int \kappa(s)ds\right]ds + n(\xi - s)\int \kappa(s)ds$$
$$+ \zeta\frac{n^2}{m^2}\left(\sinh[A]\int \kappa(s)\cosh[A]ds - \cosh[A]\int \kappa(s)\sinh[A]ds\right).$$

**Theorem 3.28.** *The position vector $\beta(s) = (\beta_1, \beta_2, \beta_3)$ of the SHCA-timelike helix of type 3 with spacelike axis $e_3$ and constant spacelike angle $\varphi$, between the unit tangent of SHCA-timelike helix and $e_3$, is obtained as*

$$\beta_1(s) = \frac{n}{m}\int\left[\int \kappa(s)\cosh[A]ds\right]ds + \frac{n}{m}(\omega - s)\int \kappa(s)\cosh[A]ds,$$
$$\beta_2(s) = \frac{n}{m}\int\left[\int \kappa(s)\sinh[A]ds\right]ds + \frac{n}{m}(\omega - s)\int \kappa(s)\sinh[A]ds,$$
$$\beta_3(s) = n\int\left[\int \kappa(s)ds\right]ds + n(\omega - s)\int \kappa(s)ds.$$

For $n = 2\sqrt{3}/3$ and $m = 2$; the graphs of spacelike slant helix $\alpha$ with timelike principal normal $N$ and curvatures $\kappa(t) = 1$, $\tau(t) = \cot\left(\frac{2\sqrt{3}}{3}t\right)$ and SHCA-timelike helix of type 1 are given in left side and right side of Fig. 2, respectively. The graphs of SHCA-timelike helices of type 2 and type 3 are given in left side and right side of Fig. 3, respectively.

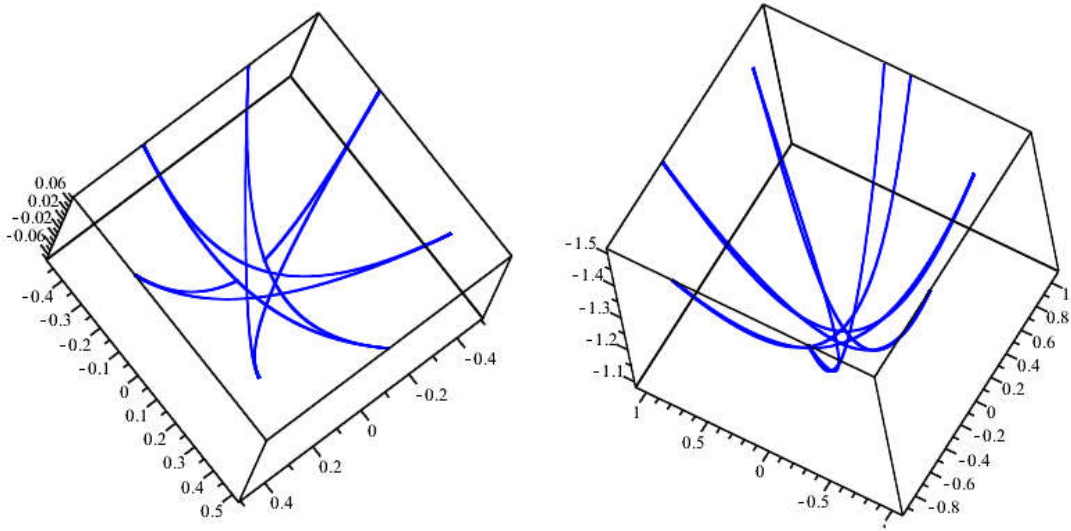

**Fig. 2.** Spacelike slant helix $\alpha$ (left) and SHCA-timelike helix of type 1 (right).

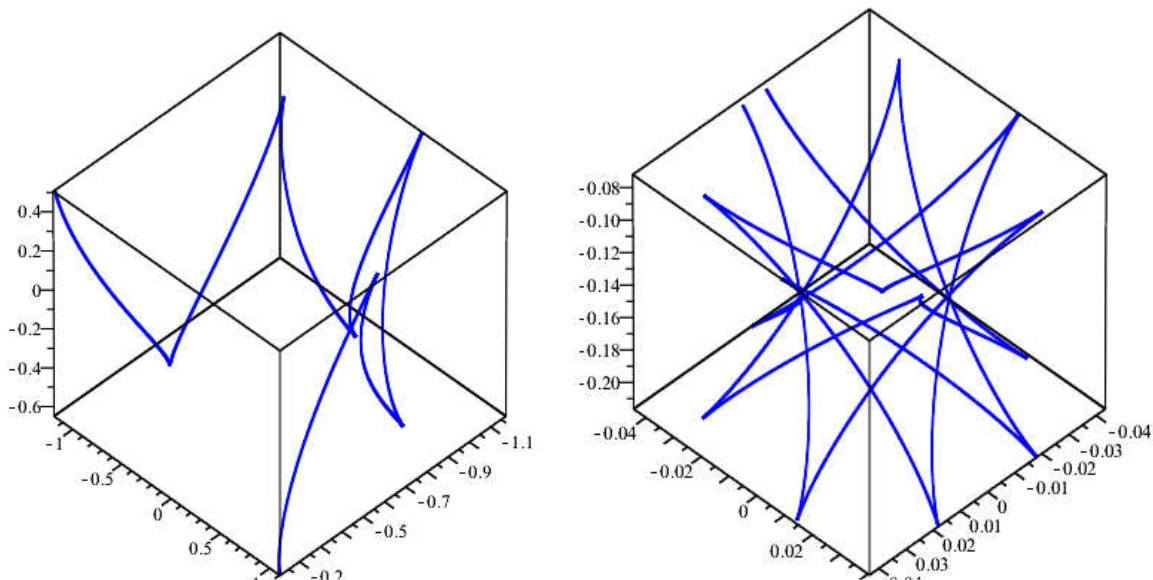

**Fig. 3.** SHCA-timelike helices of type 2(left) and type 3(right).

Moreover, if $\psi$ is a spacelike slant helix with timelike principal normal $N(s)$, the axis is timelike vector $e_1$ and constant angle between the timelike vectors $N(s)$ and $e_1$ is Lorentzian timelike angle, then the natural representation form of $\psi$ is

$$\psi(s) = n\int\left(\int\kappa(s)ds, \frac{1}{m}\int\kappa(s)\cos[A]ds, \frac{1}{m}\int\kappa(s)\sin[A]ds,\right)ds,$$

where $A = \frac{1}{n}\arcsin\left(m\int\kappa(s)ds\right)$, $\tau(s) = \pm\dfrac{m\kappa(s)\int\kappa(s)ds}{\sqrt{m^2\left(\int\kappa(s)ds\right)^2 - 1}}$ $\theta = \int\kappa(s)ds$, $m = n/\sqrt{n^2-1}$,

$n = \cosh(\varphi)$ and $\varphi$ is the constant Lorentzian timelike angle between the fixed timelike straight line $e_3$ (axis of timelike slant helix) and the tangent vector of the curve $\psi$ [7]. Then, using (26), (29) and (32), the position vectors of SHCA-timelike helices of types 1, 2 and 3 can be written similarly as given in Theorems 3.26, 3.27 and 3.28.

**Conclusions**
A new type of non-lightlike general helix called associated non-lightlike helix are introduced. These are special helices related to a non-lightlike reference curve and have different types according to reference curve. If the main curve is a timelike helix or a spacelike helix with a timelike binormal(respectively, spacelike slant helix with timelike principal normal), then associated timelike helix is called helix-connected associated (HCA)(respectively, slant helix-connected associated (SHCA)) timelike helix. Furthermore, some special types of HCA- and SHCA-timelike helices are defined and position vectors of these curves are introduced. The similar definition can be given for associated spacelike helices and corresponding theorems can be introduced by considering Lorentzian characters of associated curves.